\documentclass[12pt]{article}
\usepackage{amssymb,latexsym}
\def\ZZ{{\mathbb Z}}
\def\QQ{{\mathbb Q}}
\def\RR{{\mathbb R}}
\def\CC{{\mathbb C}}

\newtheorem{formula}{}[section]
\newtheorem{proposition}[formula]{Proposition}
\newtheorem{definition}[formula]{\indent Definition}
\newtheorem{corollary}[formula]{\indent Corollary}
\newtheorem{remark}[formula]{\indent Remark}
\newtheorem{lemma}[formula]{\indent Lemma}
\newtheorem{theorem}[formula]{\indent Theorem}

\newtheorem{example}[formula]{Example}

\def\thrm{\begin{theorem}}
\def\thrml#1{\begin{theorem}\label{#1}}
\def\ethrm{\end{theorem}}
\def\rmrk{\begin{remark}}
\def\rmrkl#1{\begin{remark}\label{#1}}
\def\ermrk{\end{remark}}
\def\dfntn{\begin{definition}}
\def\dfntnl#1{\begin{definition}\label{#1}}
\def\edfntn{\end{definition}}
\def\nmrt{\begin{enumerate}}
\def\enmrt{\end{enumerate}}

\def\qtn{\begin{equation}}
\def\qtnl#1{\begin{equation}\label{#1}}
\def\eqtn{\end{equation}}
\def\lmm{\begin{lemma}}
\def\lmml#1{\begin{lemma}\label{#1}}
\def\elmm{\end{lemma}}
\def\crllr{\begin{corollary}}
\def\crllrl#1{\begin{corollary}\label{#1}}
\def\ecrllr{\end{corollary}}

\begin{document}
\title{}
\date{}
\maketitle
\vspace{-0,1cm} \centerline{\bf Tropical Newton-Puiseux
polynomials}
\vspace{7mm}
\author{
\centerline{Dima Grigoriev}
\vspace{3mm}
\centerline{CNRS, Math\'ematique, Universit\'e de Lille, Villeneuve
d'Ascq, 59655, France} \vspace{1mm} \centerline{e-mail:\
dmitry.grigoryev@math.univ-lille1.fr } \vspace{1mm}
\centerline{URL:\ http://en.wikipedia.org/wiki/Dima\_Grigoriev} }

\begin{abstract}
We introduce tropical Newton-Puiseux polynomials admitting rational
exponents. A resolution of a tropical hypersurface is defined by
means of a tropical Newton-Puiseux polynomial. A polynomial
complexity algorithm for resolubility of a tropical curve is
designed. The complexity of resolubility of tropical prevarieties of
arbitrary codimensions is studied. Tropical Newton-Puiseux rational functions are introduced, and we prove that any tropical polynomial has a resolution in tropical 
Newton-Puiseux rational functions (this can be treated as a tropical analog of the algebraic closedness of the field of Newton-Puiseux series).
\end{abstract}

\section*{Introduction}

Recall (see e.~g. \cite{MS2015tropical}) that in the tropical semiring
$\oplus$ denotes $\min$ and $\otimes$ denotes the (classical) addition $+$.
As examples of tropical semirings one can take $\ZZ,\, \RR$. A {\it tropical}
(respectively, {\it tropical Laurent) monomial} has the form
$$a\otimes x^{\otimes I}:= a\otimes x_1^{\otimes i_1}\otimes \cdots \otimes
x_n^{\otimes i_n}$$
\noindent where $a\in \RR$ and $0\le i_1,\dots,i_n\in \ZZ$ (respectively,
$i_1,\dots,i_n\in \ZZ$). Thus, classically $a\otimes x^{\otimes I}$ equals a
linear function $a+\sum_{1\le j\le n} i_j\cdot x_j$. A {\it tropical polynomial}
 $f$
has the form $\bigoplus _I a_I \otimes x^{\otimes I}$, being classically a
convex piece-wise linear function.

A vector ${\bf x}=(x_1,\dots,x_n)\in \RR^n$ is a {\it tropical root} of $f$ if the
minimum of $a_I \otimes {\bf x}^{\otimes I}$ is attained at least for two
different tropical monomials of $f$. The set of all tropical roots of $f$
constitute a {\it tropical hypersurface} $T(f)\subset \RR^n$ being a finite
union of polyhedra of dimensions $n-1$.

We extend the concept of a tropical polynomial by allowing the exponents
$i_1,\dots, i_n$ to be rational calling it a {\it tropical Newton-Puiseux
polynomial}. Assume that
\begin{eqnarray}\label{5}
f=\bigoplus _{0\le i\le d} f_i\otimes y^{\otimes i}
\end{eqnarray}
for some tropical polynomials $f_0,\dots,f_d$ in the variables
$x_1,\dots,x_n$. We call a Newton-Puiseux polynomial $y$ a {\it
resolution} of $f$ (or of the tropical hypersurface $T(f)$) if for
any point ${\bf x}\in \RR^n$ the point $({\bf x},y({\bf x}))\in \RR^{n+1}$
provides a tropical root of $f$ (the formal definitions one can find
below in Section~\ref{resolution}).

This resembles Newton-Puiseux series from algebraic geometry with the difference
 that we consider finite supports since in the tropical semiring one takes
$\min$. Tropical Newton-Puiseux polynomials can be viewed as a tropical analog of algebraic functions.

In Section~\ref{resolution} we show that the set of all the
resolutions of a tropical hypersurface is finite and closed under
taking $\min$. Thus, there exists a minimal resolution, and in case
of a monic tropical polynomial
$$f=y^{\otimes d}\oplus \bigoplus _{0\le i < d} f_i \otimes y^{\otimes i}$$
\noindent we provide a simple formula for the minimal resolution. In
addition, a geometric description of resolutions is given. Also we
show that the resolubility of a tropical hypersurface  belongs to
the complexity class $NP$.

In Section~\ref{divisibility} a polynomial (bit-size) complexity algorithm is exhibited
for resolving degree 1 tropical polynomials of the form $f_1\otimes y \oplus f_0$, which is equivalent to the divisibility of $f_0$ by $f_1$.

In Section~\ref{algorithm} we design a polynomial (bit-size) complexity algorithm for testing resolubility of a tropical curve in a real space of a fixed dimension, moreover the algorithm provides
a succinct description of the set of all the resolutions.

In Section~\ref{several} we study the problem of resolubility of a system of tropical polynomials in a single variable $x$ and in several indeterminates
$y_1,\dots, y_s$ and establish its $NP$-hardness.

In Section~\ref{rational} we study tropical Newton-Puiseux rational
functions, being tropical quotients (or in other words, the
classical subtraction) of pairs of tropical Newton-Puiseux
polynomials. We prove (see Remark~\ref{non-monic}) a tropical analog of the algebraic closedness of the field of Newton-Puiseux series, namely, that any tropical polynomial $f$ (see (\ref{1})) has a resolution in tropical Newton-Puiseux rational functions, and moreover, provide an explicit formula
$$\bigoplus _{0\le i\le d} (f_{d-i}\oslash f_d)^{\otimes (1/i)}$$ 
\noindent for the minimal resolution of $f$. Also an algorithm is suggested which tests resolubility of a
tropical curve by means of tropical Newton-Puiseux rational
functions. The complexity of the algorithm is polynomial for a fixed
dimension of the ambient space.

\section{Resolution of a tropical hypersurface}\label{resolution}

Let an algebraic (classical) equation
\begin{eqnarray}\label{1}
F:=\sum_{0\le i\le d} F_i\cdot Y^i=0
\end{eqnarray}
where the coefficients $F_i\in K[X_1,\dots,X_n]$ for the field
$K=\CC((t^{1/\infty}))$ of Newton-Puiseux series, have a Laurent
polynomial solution
\begin{eqnarray}\label{2}
Y=\sum_I A_I\cdot X^I
\end{eqnarray}
with a finite sum over multiindices $I\in \ZZ^n$ and the
coefficients $A_I\in K$.

Denote the {\it tropicalization}
\begin{eqnarray}\label{3}
Trop(Y):=\bigoplus _I Trop(A_I) \otimes X^{\otimes I}
\end{eqnarray}
where for a Newton-Puiseux series $A_I=\sum _{0\le j<\infty}
b_j\cdot t^{s_j/q}$ with $b_j\in \CC,\, b_0\neq 0$ and increasing
integers $s_0<s_1<\dots$ its tropicalization $Trop(A_I):=s_0/q\in
\QQ$.

\begin{remark}\label{motivation}
$Trop(Y)$ is a solution of the tropical equation
\begin{eqnarray}\label{4}
\bigoplus _{0\le i\le d} Trop(F_i)\otimes (Trop(Y))^{\otimes i}
\end{eqnarray}
This means that for any point ${\bf x}=(x_1,\dots,x_n)\in \RR^n$ the
minimal value of $Trop(F_i)\otimes (Trop(Y))^{\otimes i}$ at ${\bf x}$
for $0\le i\le d$ is attained at least for two different $0\le
i_1<i_2\le d$.
\end{remark}

\begin{remark}\label{integer}
Observe that the validity of (\ref{4}) would not change if one
multiplies all the rational coefficients in $Trop(F_i),\, 0\le i\le
d$ by their common denominator $m$ and simultaneously all
$Trop(A_I)$ (see (\ref{3})) by $m$ to make all the coefficients in
$Trop(F_i),\, 0\le i\le d$ integers.
\end{remark}

Remark~\ref{motivation} motivates the following definition.

\begin{definition}\label{definition}
A tropical hypersurface $T(f)\subset \RR^{n+1}$ defined by a
tropical polynomial (\ref{5}) where $f_i$ are tropical polynomials
in the variables $x_1,\dots,x_n$ with integer coefficients (cf.
Remark~\ref{integer}) has a {\it resolution} being a tropical
Newton-Puiseux polynomial
\begin{eqnarray}\label{6}
y=\bigoplus _I a_I \otimes x^{\otimes I}
\end{eqnarray}
for a finite sum over multiindices $I\in \QQ^n$ and $a_I\in \QQ$, if
for any point ${\bf x}=(x_1,\dots,x_n)\in \RR^n$ the minimal value among
$f_i\otimes y^{\otimes i},\, 0\le i\le d$ (treated as piece-wise
linear functions) at ${\bf x}$ is attained at least for two different
$0\le i_1<i_2\le d$.
\end{definition}

Denote by $N$ the common denominator of all the rational coordinates of multiindices $I$ from (\ref{6}). Then $y^{\otimes N}$ is a tropical (Laurent) polynomial
which equals classically to $N\cdot \min_I\{a_I+i_1x_1+\cdots +i_nx_n\}$. 

\begin{proposition}
Let $y$ be a resolution of $f$ (see (\ref{5}), (\ref{6})), then
$({\bf x},y({\bf x}))\in T(f)$.
\end{proposition}

\begin{example}
The tropical polynomial $f=y\oplus x\oplus 0$ has a resolution
$y=x\oplus 0$. Its graph $\{({\bf x},y({\bf x})): {\bf x}\in \RR\} \subset T(f)\subset
\RR^2$ consists of two half-lines, while the tropical curve $T(f)$
consists of three half-lines.
\end{example}

\begin{proposition}\label{minimum}
Let $y$ (see (\ref{6})) and $\bigoplus _I b_I \otimes x^{\otimes I}$
be resolutions of (\ref{5}). Then $\bigoplus _I (a_I \oplus
b_I)\otimes x^{\otimes I}$ is also a resolution of (\ref{5}).
\end{proposition}

The proof follows from an observation that for any point ${\bf x}\in
\RR^n$ the minimum on the tropical monomials after opening the
parenthesis in a power $y^{\otimes i}$ (see (\ref{6})) is attained
on the powers of the kind $(a_I\otimes {\bf x}^{\otimes I})^{\otimes
i}$. $\Box$ \vspace{2mm}

Below in Remark~\ref{np}
we show that there is at most a finite number of resolutions of
(\ref{5}). Hence according to Proposition~\ref{minimum}, there
exists a minimal resolution.

\begin{proposition}\label{monic}
If $$f=y^{\otimes d} \oplus \bigoplus _{0\le i<d} f_i \otimes y^{\otimes i}$$
\noindent (see (\ref{5})) is monic then $$y=\bigoplus _{1\le i \le
d} f_{d-i} ^{\otimes (1/i)}$$ \noindent is the minimal resolution.
\end{proposition}

{\bf Proof}. For any point ${\bf x}\in \RR^n$ the minimal $y_0\in \RR$
such that $({\bf x},y_0)$ belongs to the tropical hypersurface
$T(f)\subset \RR^{n+1}$ satisfies a (classical) equation $d\cdot
y_0=f_i({\bf x})+i\cdot y_0$ for suitable $0\le i < d$ (due to analyzing
the Newton polygon). $\Box$ \vspace{2mm}

Note also that if $f_{d-i}=\bigoplus _J c_J \otimes x^{\otimes J}$
then $f_{d-i}^{\otimes (1/i)} = \bigoplus _J (c_J/i)\otimes
x^{\otimes (J/i)}$.

\begin{remark}\label{non-monic}
When $f$ is not monic, a resolution does not necessary exist as in
the example $f=(x\oplus 0)\otimes y \oplus 0$. On the other hand,
one can write a similar formula
$$y=\bigoplus _{1\le i\le d} (f_{d-i} \oslash f_d)^{\otimes (1/i)}$$
\noindent where $\oslash$ stays for the tropical division, i.~e. the
classical subtraction. In this case $y$ is not necessary a convex
function, while being  piece-wise linear (we call them tropical Newton-Puiseux rational functions, see
Section~\ref{rational}), and $y$ provides the minimal resolution of $f$. This can be treated as a tropical analog of the algebraic closedness of the field of Newton-Puiseux
series.
\end{remark}

Now we proceed to a geometric description of resolutions. Let
(\ref{6}) be a resolution of (\ref{5}). Assume that for some $I$ the
(convex) polyhedron $M_I\subset \RR^n$ of points at which the
(tropical) monomials $\{a_J\otimes x^{\otimes J}\}_J$ of $y$ attain
the minimum for $a_I\otimes x^{\otimes I}$, has the full dimension
$n$. Observe that if $M_I$ has a dimension less than $n$ one can
discard the monomial $a_I\otimes x^{\otimes I}$ from $y$.

Assume that for some $0\le i_1<i_2\le d$ and a pair of monomials
$c_{i_1,I_1}\otimes x^{\otimes I_1},\, c_{i_2,I_2}\otimes x^{\otimes
I_2}$ from the polynomials $f_{i_1},\, f_{i_2}$, respectively, it
holds
\begin{eqnarray}\label{7}
I_1+i_1\cdot I=I_2+i_2\cdot I;\, c_{i_1,I_1}+i_1\cdot
a_I=c_{i_2,I_2}+i_2\cdot a_I,
\end{eqnarray}
in other words, the monomials
$$(c_{i_1,I_1}\otimes x^{\otimes I_1})\otimes (a_I\otimes x^{\otimes
I})^{\otimes i_1}= (c_{i_2,I_2}\otimes x^{\otimes I_2})\otimes
(a_I\otimes x^{\otimes I})^{\otimes i_2}$$ \noindent coincide.
Consider the convex polyhedron $M_{I,i_1,I_1,i_2,I_2} \subset M_I$
of the points from $M_I$ at which the minimum of the monomials
$(c_{i,I}\otimes x^{\otimes I})\otimes (a_I\otimes x^{\otimes
I})^{\otimes i}$ for the monomials $c_{i,I}\otimes x^{\otimes I}$
from $f_i,\, 0\le i\le d$ is attained for $(c_{i_1,I_1}\otimes
x^{\otimes I_1})\otimes (a_I\otimes x^{\otimes I})^{\otimes i_1}$.
We get the following lemma.

\begin{lemma}\label{lemma}
Let (\ref{6}) be a resolution of (\ref{5}) and the polyhedron
$M_I\subset \RR^n$ have the full dimension $n$. Then the polyhedra
$M_{I,i_1,I_1,i_2,I_2}$ having the full dimension $n$ constitute a
partition of $M_I$, i.~e. every two elements of the partition either
coincide or intersect by a set (face) of dimension less than $n$.
\end{lemma}

It would be interesting to clarify, how many resolutions a
tropical hypersurface might have?

Let the tropical degrees $trdeg (f_i)\le D,\, 0\le i\le d$.

\begin{remark}\label{np}
The problem of
resolving a tropical polynomial (\ref{5}) belongs to the complexity
class $NP$. This follows from the observation that each coefficient
$a_I$ satisfies (\ref{7}) (or equals infinity), and therefore, there
are at most $d^2\cdot {D+n \choose n}$ possibilities for $a_I$,
taking into account that $D+n \choose n$ bounds the number of
monomials in each $f_i$.

Note that when $f_i,\, 0\le i\le d$ are in sparse encoding, in the
latter bound one can replace $D+n \choose n$ by the number of
monomials in $f_i,\, 0\le i\le d$. Thus, the problem of resolubility
of (\ref{5}) belongs to $NP$ for both dense and sparse encodings of
(\ref{5}).
\end{remark}

It would be interesting to say more about the complexity of
resolubility of (\ref{5}).

\begin{remark}
One can extend the results of this section to an input tropical
Newton-Puiseux polynomials in place of (\ref{5}).
\end{remark}

\section{Polynomial complexity testing divisibility of tropical
polynomials}\label{divisibility}

If (\ref{5}) has degree 1, i.~e. $f=f_1\otimes y\oplus f_0$ then
according to (\ref{7}) a resolution (\ref{6}) is equivalent to the
divisibility $f_1\otimes y = f_0$ with $y$ being a tropical Laurent
polynomial. We agree that two tropical (Laurent) polynomials are
equal if they are equal as (convex piece-wise linear) functions.

We expose an algorithm for testing divisibility within polynomial
complexity. First the algorithm deletes from $f_0$ all the monomials of the form
$b\otimes x_1^{b_1}\otimes \cdots \otimes x_n^{b_n}$ which do not
change $f_0$ as a function. Geometrically, it means that the
hyperplane defined as the graph
$$\{(x_1,\dots,x_n,\sum_{1\le j\le n}b_j\cdot x_j+b) : (x_1,\dots,x_n)\in
\RR^n\}$$ \noindent of this monomial in $\RR^{n+1}$ is higher (with
respect to the last coordinate) than the polyhedron $P$ defined by
the other monomials of $f_0$ (observe that $P$ is the graph of $f_0$
as a function). The latter is a problem of linear programming. Thus,
one can suppose $f_0$ to be reduced, i.~e. do not contain
unnecessary monomials. Also we suppose that $f_1$ is reduced.

For every candidate $I=(i_1,\dots,i_n)\in \ZZ^n,\, \sum _{1\le j\le
n} |i_j| \le D$ (see (\ref{7})) to be in the support of a resolution
$y$ the algorithm calculates (again involving linear programming)
the minimal $a_I$ such that for each monomial $c\otimes x^{\otimes
C}$ of $f_1$ the hyperplane in $\RR^{n+1}$ defined by the monomial
$(c\otimes x^{\otimes C})\otimes (a_I\otimes x^{\otimes I})$ is
(non-strictly) higher than $P$.

Then $y=\bigoplus _I a_I \otimes x^{\otimes I}$ is a resolution of
$f_1\otimes y \oplus f_0$ iff for each monomial $b\otimes x^{\otimes
B}$ of $f_0$ there exists $I$ and a monomial $c\otimes x^{\otimes
C}$ of $f_1$ such that $(a_I\otimes x^{\otimes I})\otimes (c\otimes
x^{\otimes C})=b\otimes x^{\otimes B}$. Reducing further $y$ as
described above, we conclude that there is a unique reduced
resolution $y$ (provided that it does exist).

Summarizing, we obtain the following proposition.

\begin{proposition}
One can test resolubility of degree 1 tropical polynomial
$f_1\otimes y \oplus f_0$ (or equivalently, the divisibility
$f_1\otimes y =f_0$) within polynomial complexity. In case of the
divisibility the algorithm yields the unique reduced resolution $y$.
\end{proposition}

\section{Polynomial complexity algorithm for resolving tropical
curves}\label{algorithm}

Let a system of tropical polynomials
\begin{eqnarray}\label{13}
f_i,\, 1\le i\le k
\end{eqnarray}
in $n$ variables $x,y_1,\dots,y_{n-1}$ with integer coefficients
determine a tropical prevariety $T:=T(f_1,\dots,f_k)\subset \RR^n$.
Let the tropical degrees $trdeg(f_i)\le d,\, 1\le i\le k$ and the
bit-sizes of the coefficients of $f_i,\,1\le i\le k$ do not exceed
$L$.

First, the algorithm constructs $T$ as a union of polyhedra (see
e.~g. \cite{MS2015tropical}). Each of these polyhedra (including
faces of all the dimensions) is defined by specifying the monomials
of $f_i,\, 1\le i\le k$ (treated as linear functions) on which the
minima are attained (cf. e.~g. \cite{GP}). The algorithm can find
the partition of $\RR^n$ into polyhedra defined by given feasible
tuples of signs (i.~e. either the positive, either the negative or
zero) of all the differences of the monomials of $f_i,\ 1\le i\le k$
(in other words, by all the feasible orderings of the monomials of
$f_i,\, 1\le i\le k$). Namely, the algorithm finds the partition by
recursion on the number of the differences. If for a current subset
of the differences the partition of $\RR^n$ w.r.t. this subset is
already constructed, the algorithm picks up the next difference and
for each element (being a polyhedron) of the current partition
verifies which signs of the picked up difference are feasible on
this polyhedron (with the help of linear programming). Thereupon,
the algorithm discards the unfeasible tuples of signs, which
completes the recursive step.

The number of the elements of a current partition at every step of
the recursion is bounded by $$n^2\cdot 2^n\cdot {{k\cdot {d+n
\choose n}^2} \choose n} < k^n\cdot d^{2\cdot n^2}$$ \noindent due
to the Buck's formula on the number of faces in an arrangement of
hyperplanes \cite{fukuda}. Hence the complexity of the recursion is
bounded by a polynomial in $L,\, k^n,\, d^{n^2}$ because the
algorithm invokes linear programming the number of times being
polynomial in $k^n\cdot d^{n^2}$.

Since the tropical prevariety $T$ is a union of appropriate subset
of the elements of the constructed partition of $\RR^n$, we get the
following proposition.

\begin{proposition}\label{prevariety}
There is an algorithm which constructs the tropical prevariety
$T(f_1,\dots,f_k)\subset \RR^n$ determined by (\ref{13}) within the
complexity polynomial in $L,\, k^n,\, d^{n^2}$.
\end{proposition}

Now we assume that $\dim T =1$, thus $T$ is a tropical curve. We
design an algorithm which verifies the resolubility of $T$, i.~e.
whether there exist tropical Newton-Puiseux polynomials
$y_1(x),\dots,y_{n-1}(x)$ assuring a resolution of (\ref{13}). The
latter is equivalent to that every piece-wise linear function
$y_j,\, 1\le j\le n-1$ is convex.

The algorithm produces a directed graph $G$ whose vertices being the
edges of $T$ (including the unbounded ones) not lying in a
hyperplane of the form $x=c$. Two edges $e^{(-)},\, e^{(+)}$ of $T$
(being vertices of $G$) with the same endpoint of the kind
$$e^{(-)}=((x^{(-)},y_1^{(-)},\dots,y_{n-1}^{(-)}),\,
(x,y_1,\dots,y_{n-1})),$$
$$e^{(+)}=((x,y_1,\dots,y_{n-1}),\,(x^{(+)},y_1^{(+)},\dots,y_{n-1}^{(+)}))$$
\noindent are linked by an edge directed from $e^{(-)}$ to $e^{(+)}$
in $G$ if
\begin{eqnarray}\label{14}
x^{(-)}<x<x^{(+)};\quad \frac{y_j-y_j^{(-)}}{x-x^{(-)}} \ge
\frac{y_j-y_j^{(+)}}{x-x^{(+)}},\, 1\le j\le n-1.
\end{eqnarray}
When $e^{(-)}$ (respectively, $e^{(+)}$) is unbounded with an
endpoint $(x,y_1,\dots,y_{n-1})$ (so, is a half-line), which we call
unbounded from the left, we take an arbitrary point of $e^{(-)}$
with $x^{(-)}<x$ (respectively, if $e^{(+)}$ is a half-line, we take
a point of $e^{(+)}$ with $x^{(+)}>x$, and we call $e^{(+)}$
unbounded from the right). When an edge of $T$ has no endpoints, so
is a line, it provides a resolution of $T$.

After that the algorithm produces a subset $S$ of the vertices of
$G$. It starts with including into $S$ all the edges of $T$ (so, the
vertices of $G$) unbounded from the left (denote this set by $S_0$).
Thereupon, the algorithm includes into $S$ all the vertices of $G$
reachable from $S_0$. If a vertex of $G$ corresponding to an edge of
$T$ unbounded from the right, belongs to $S$, a path in $G$ leading
to such a vertex from $S_0$ provides a resolution of $T$ (i.~e. each
piece-wise linear function $y_j(x),\, 1\le j\le n-1$ corresponding
to the path, is convex due to (\ref{14})). Moreover, the paths in
$G$ from $S_0$ to the vertices corresponding to the edges of $T$
unbounded from the right, are in a bijective correspondence with the
resolutions of $T$.

Summarizing and taking into account Proposition~\ref{prevariety}, we
obtain the following theorem.

\begin{theorem}
There is an algorithm which tests resolubility of a tropical curve
$T\subset \RR^n$ determined by (\ref{13}), and in case of the
resolubility yields a resolution. The complexity of the algorithm is
polynomial in $L,\,k^n,\, d^{n^2}$. In particular, the complexity is
polynomial for a fixed ambient dimension $n$.
\end{theorem}

\begin{remark}
Let a system of tropical polynomials of the form (\ref{13}) depend
on the variables $x_1,\dots,x_m,y_1,\dots,y_n$ and the tropical
prevariety $T\subset \RR^{m+n}$ have dimension $m$. Then one can try
different subsets of all $m$-dimensional faces of $T$ as candidates
to constitute a graph of a resolution $$(x_1,\dots,x_m)\to
(x_1,\dots,x_m,y_1(x_1,\dots,x_m),\dots,y_n(x_1,\dots,x_m)) \in T$$
\noindent of $T$ similar to Remark~\ref{np}. The latter, in fact,
means that firstly, the projections of the chosen $m$-dimensional
faces on $\RR^m$ with the coordinates $x_1,\dots,x_m$ form a
partition of $\RR^m$ and secondly, that each piece-wise linear
function $y_j(x_1,\dots,x_m),\, 1\le j\le n-1$ is convex.
\end{remark}

\section{Resolution of systems of tropical polynomials with several
indeterminates}\label{several}

In this section we consider systems of tropical polynomials (instead
of a single polynomial (\ref{5})) in one variable $x$ and several
indeterminates $y_1,\dots,y_s$. Thus, in a resolution (cf.
(\ref{6})) each $y_i$ is a tropical Newton-Puiseux polynomial. We
show the following proposition.

\begin{proposition}
The problem of resolubility of a system of tropical polynomials in a
single variable and in several indeterminates is $NP$-hard.
\end{proposition}

{\bf Proof}. We reduce 3-SAT to the problem under consideration, so
we construct a system $R$ of tropical polynomials. For an instance
of 3-SAT problem in $n$ variables $u_1,\dots,u_n$ we introduce
indeterminates $y_1,\dots,y_n,\, z_1,\dots,z_n$ and add to $R$
tropical polynomials
\begin{eqnarray}\label{8}
y_i\otimes z_i \oplus x,\, 1\le i\le n
\end{eqnarray}
Formula (\ref{8}) means that the resolutions of $y_i$ and of $z_i$
are both monomials in $x$. Informally, $0=x^{\otimes 0}$ encodes the
truth and $x=x^{\otimes 1}$ encodes the falsity, $y_i$ corresponds to
$u_i$ and $z_i$ corresponds to $\neg u_i$.

For every $j$-th 3-clause of the 3-SAT formula, say, $u_m \vee \neg
u_k \vee u_l$ we add to $R$ the following tropical (linear)
polynomials
\begin{eqnarray}\label{9}
y_m\oplus z_k\oplus y_l\oplus v_j;
\end{eqnarray}
\begin{eqnarray}\label{10}
v_j\oplus x^{\otimes 1}\oplus w_j;
\end{eqnarray}
\begin{eqnarray}\label{11}
w_j\oplus x^{\otimes 1} \oplus 0
\end{eqnarray}
with indeterminates $v_j,\, w_j$. Note that (\ref{11}) ensures that
in a resolution the reduced $w_j=x^{\otimes 1} \oplus 0$, then
(\ref{10}) ensures that the reduced $v_j$ contains the constant
monomial $0$ (and possibly, monomials of the form $c\otimes
x^{\otimes b}$ with $0<b\le 1,\, c\ge 0$). Finally, (\ref{9})
ensures that one of the resolutions of $y_m,\, z_k,\, y_l$ equals
$0$.

Thus, existence of a resolution of the system $R$ for all $j$
implies the solvability of the initial 3-SAT formula.

The converse is obvious: for a Boolean vector $(u_1,\dots,u_n)$
providing a solution of the initial 3-SAT formula put $y_i=0,\,
z_i=x^{\otimes 1}$ when $u_i$ is true and $y_i=x^{\otimes 1}, \,
z_i=0$ for $u_i$ being false. Thereupon put $v_j=y_m\oplus z_k\oplus
y_l$. $\Box$

We mention that the problem of solvability of a system of tropical
polynomials is $NP$-complete \cite{theobald06frontiers}.

It would be interesting to understand more about the complexity of
the problem under consideration in this section.
 \vspace{2mm}

 \section{Tropical Newton-Puiseux rational
 functions}\label{rational}

 Any tropical Newton-Puiseux rational function $f_1\oslash f_2$ where
 $f_1,\, f_2$ are tropical Newton-Puiseux polynomials, is a
 piece-wise linear (continuous) function (cf.
 Remark~\ref{non-monic}). The converse is also true (see e.~g.
 \cite{bittner}, \cite{ks}): any piece-wise linear continuous
 function is a difference of two piece-wise linear convex functions.
 In \cite{ovchinnikov} an algorithm is suggested which represents a
 piece-wise linear function as a difference of piece-wise linear
 convex functions with the complexity bound being exponential. In
 case of one-variable functions a polynomial complexity algorithm
 for this problem is exhibited in \cite{dobkin}.




 Let a tropical curve $T\subset \RR^n$ be determined by a system
 (\ref{13}). As in Section~\ref{algorithm} the algorithm finds $T$.
 Thereupon, similar to Section~\ref{algorithm} constructs a graph
 which comprises all the paths consisting of the edges of $T$ of the
 form $$\{(x_l,y_1^{(l)},\dots,y_{n-1}^{(l)}),\,
 (x_{l+1},y_1^{(l+1)},\dots,y_{n-1}^{(l+1)})\, : \, 0\le l\le s\}$$
 \noindent where $x_0:=-\infty< x_1 < \cdots < x_s <
 x_{s+1}:=\infty$, thus, this path contains $s+1$ edges. The
 difference with Section~\ref{algorithm} is that now we do not
 impose a requirement on convexity.

 The algorithm can pick up any such path (provided that it does
 exist), then this path yields $n-1$ piece-wise linear functions $y_i(x),\, 1\le i< n$. Making use of \cite{dobkin} the algorithm
 represents $y_i(x)=g_i(x)-h_i(x)$ with piece-wise linear convex
 functions $g_i,\, h_i$. This produces a tropical Newton-Puiseux
 rational function resolution of $T$.

 Summarizing and invoking the complexity bounds from
 Section~\ref{algorithm}, we get the following proposition.

 \begin{proposition}
There is an algorithm which tests resolubility of a tropical curve
determined by (\ref{13}) by means of tropical Newton-Puiseux
rational functions within the complexity polynomial in $L,\,k^n,\,
d^{n^2}$. The algorithm yields a resolution, provided that it does
exist. Therefore, the complexity is polynomial for a fixed dimension
$n$ of the ambient space.
 \end{proposition}

 It would be interesting to estimate the complexity of resolubility
 of tropical prevarieties or arbitrary dimensions by means of tropical
 Newton-Puiseux rational functions. \vspace{2mm}

{\bf Acknowledgements}. The author is grateful to the grant RSF 16-11-10075
and to MCCME for inspiring atmosphere.

\end{document}